# INVARIANT *P*-VALUES FOR MODEL CHECKING

By Michael Evans and Gun Ho Jang

*University of Toronto*

*P*-values have been the focus of considerable criticism based on various considerations. Still, the *P*-value represents one of the most commonly used statistical tools. When assessing the suitability of a single hypothesized distribution, it is not clear that there is a better choice for a measure of surprise. This paper is concerned with the definition of appropriate model-based *P*-values for model checking.

**1. Introduction.** The use of *P*-values is common in statistical practice. Despite this, it is reasonable to say that the logical foundations for the *P*-value are somewhat weak. This has lead to a variety of criticisms of *P*-values and even to doubts as to their correctness; see, for example, the discussions in [3–5, 9, 12, 13] and [15]. For example, [9] is concerned with the use of the 5% cut-off as a standard for determining whether or not a result is "significant," while [3] and [4] argue that frequentist *P*-values can be misleading.

While various arguments have been advanced for alternatives to *P*-values, there are situations where the use of *P*-values seems unavoidable. One such context arises with model checking, where we have observed data $x_0 \in \mathcal{X}$ and want to assess whether or not $x_0$ is a plausible value from a fixed probability measure $P$. For example, in model checking, $P$ could arise as the conditional distribution of the observed data given a minimal sufficient statistic or as the distribution of an ancillary statistic such as a function of residuals. If the *P*-value leads us to doubt that $x_0$ could have arisen from $P$, then we also have reason to doubt the underlying statistical model.

So, the basic problem we consider is to determine a measure of how surprising the observed value $x_0$ is as a possible value from $P$. A common approach to this is to say that we need to prescribe a real-valued discrepancy statistic $T: \mathcal{X} \to R^1$, so that, in some sense, $T(x)$ measures how divergent the value $x$ is, and then to compute the $P$-value

$$(1) \qquad\qquad P(T(x) \geq T(x_0)).$$

If (1) is small, then we are lead to doubt $P$. In general, no guidance is provided as to how the statistic $T$ is to be chosen with respect to $P$. For example, we do not have a likelihood ratio available in this situation as a possible choice of $T$. Further, it can be noted that some restrictions on $T$ are necessary if (1) is to have an appropriate interpretation. In particular, the right tail of the distribution of $T$ should be the only region that has relatively low probability. Otherwise, we could have a value of $T(x_0)$ in the left tail of the induced measure $P_T$ or near a shallow anti-mode of $P_T$, that leads to a reasonable value of (1)—and yet $T(x_0)$ could still be considered as surprising.

In the case where $T$ has a discrete probability distribution, there is a natural definition of a $P$-value that avoids these problems, namely,

$$(2) \qquad\qquad P(p_T(T(x)) \leq p_T(T(x_0))),$$

where $p_T$ is the probability function of $P_T$. In this case, we see that (2) is the probability of obtaining a value of $T$ with probability of occurrence no greater than the probability of what was actually observed. If this probability is small, then $T(x_0)$ is a surprising value. We immediately see that (2) will identify values of $T$ in either tail, or near shallow anti-modes, as being surprising. In fact, there is no need to require that $T$ be real-valued for (2) to make sense.

While (2) seems like a very natural definition, a serious problem arises when we attempt to generalize this to the situation where the distribution of $T$ is absolutely continuous, say with density $f_T$ with respect to some support measure on the range space $\mathcal{T}$ of $T$. Intuitively, we would like to use the analog of (2) given by

$$(3) \qquad\qquad P(f_T(T(x)) \leq f_T(T(x_0))).$$

But now suppose that we have a 1–1, smooth transformation $W : \mathcal{T} \to \mathcal{T}$. The density of $W$ is $f_W(w) = f_T(W^{-1}(w)) J_W(W^{-1}(w))$, where $J_W(t)$ is the reciprocal of the Jacobian determinant of $W$ at $t$. The $P$-value based on the density of $W$ is $P(f_W(W \circ T(x)) \leq f_W(W \circ T(x_0))) = P(f_T(T(x)) J_W(T(x)) \leq f_T(T(x_0)) J_W(T(x)))$ and it is clear that, unless $J_W(T(x)))$ is constant, this will not equal (3). In fact, these values can be quite different. We refer to this as the noninvariance of the $P$-value given by (3).

It is the goal of this paper to provide a definition of a $P$-value for the absolutely continuous case that is invariant, one which proceeds logically from the natural definition in the discrete case as given by (2). We base our argument on the intuitively reasonable idea that any absolutely continuous model is, in fact, an approximation of an underlying model that is discrete, as it is based on a measurement process for a response that has finite accuracy. When we take this approximation into account and adjust for any volume distortions induced by transformations such as discrepancy statistics $T$, we are led to an invariant $P$-value. The central idea is that we do not want inferences in statistical problems to be based on changes in volume as induced by transformations and so we must adjust the density $f_T$ appropriately in (3) to avoid this. We provide the basic definition and argument in Section 2, along with examples that support our approach. In Section 3, we apply these results to model checking problems. In Section 4, we draw some conclusions.

**2. A general $P$-value.** To start, we consider a response $x$ taking values in an open subset $\mathcal{X}$ of $R^k$ and develop an appropriate $P$-value when $T(x) \equiv x$. We suppose that the probability measure $P$ has density $f$ with respect to support measure $\mu$ on $\mathcal{X}$.

In measure-theoretic terms, a density $f$, with respect to a support measure $\mu$, is seen simply as a device to compute probabilities. In statistical contexts, however, a density plays a somewhat more significant role. For example, if $f(x_1) > f(x_2)$, then we want to say that the probability of $x_1$ occurring is greater than the probability of $x_2$ occurring. For this to hold, we cannot allow $f$ to be defined in an arbitrary fashion. In effect, we need to have that $P(A)/\mu(A) \rightarrow f(x)$ as the set $A$ converges to $\{x\}$ since this implies that $P(A) \approx f(x)\mu(A)$ when $A$ is close to $\{x\}$. Further, to compare the probabilities of two points $x_1$ and $x_2$, we need $A_i \rightarrow \{x_i\}$ with $\mu(A_1) = \mu(A_2)$ and then, for example, we can say that the probability of $x_1$ occurring is greater than the probability of $x_2$ occurring when $f(x_1) > f(x_2)$. The mathematics of making this precise is discussed in, for example, [14], under the topic of differentiating one measure with respect to another.

It seems natural to choose $\mu = \mu_k$, where $\mu_k$ is Euclidean volume on $\mathcal{X}$, as it weights sample points equally and so $f(x)$ expresses the essence of how the probability measure is behaving at $x$. This is analogous to using counting measure as the support measure in the discrete case since $f(x)$ then has a direct interpretation as the probability of $x$.

We now consider the essential discreteness of the observational process. Suppose that this translates into a value $x$ lying in a set $B_n(x)$ such that $\{B_n(x) : x \in \mathcal{X}\}$ forms a partition of $\mathcal{X}$ with $\mu_k(B_n(x))$ finite and constant in $x$ and such that $B_n(x)$ shrinks "nicely" (see [14]) to $x$ as $n \rightarrow \infty$. We then have that $P(B_n(x))/\mu_k(B_n(x)) \rightarrow f(x)$ as $n \rightarrow \infty$ as long as $f$ is continuous at $x$. So, for $n$ large, $P(B_n(x)) \approx f(x)\mu_k(B_n(x))$ and $f(x)$ serves as surrogate for the probability of $x$, at least when we are comparing the probabilities of different values of $x$ occurring. Note that the constancy of $\mu_k(B_n(x))$ in $x$ is necessary for this interpretation of $f(x)$. As a particular example of this, suppose that $\mathcal{X} = R$, that we partition $R^1$ using $\{((i-1)/n, i/n] : i \in Z\}$ and that $B_n(x)$ is the set $((i-1)/n, i/n]$ that contains $x$.

Rather than observing $x$, the essential discreteness of the problem means that we will observe some $x_n(x) \in B_n(x)$ and the probability of observing $x_n(x)$ is $P(B_n(x))$. Note that, implicitly, $x_0$ is one of the values assumed by $x_n$. Then, for the discrete response variable $x_n$, the appropriate $P$-value (2) is given by

$$(4) \qquad \sum_{\{x_n(x)\,:\,P(B_n(x)) \le P(B_n(x_0))\}} P(B_n(x)).$$

We then want to show that (3), with $f_T = f$, approximates (4).

While such a result seems intuitively plausible, a general proof is not straightforward. We require some regularity conditions as we cannot expect such an approximation to hold if we allow $f$ and the partition $\{B_n(x) : x \in \mathcal{X}\}$ to be too general. For this, we use the theory of contented sets and functions, as discussed in [10], where it is used to develop the Riemann integral. Essentially, a bounded set $A$ is *contented* if its $\mu_k$-measure can be approximated arbitrarily closely by the $\mu_k$-measure of a finite union of disjoint rectangles contained in $A$ and also by the $\mu_k$-measure of a finite union of disjoint rectangles containing $A$. A bounded function $f$ with compact support is *contented* if it can be approximated arbitrarily closely by step functions based on contented sets. Further, we say that a function $f$ is *locally constant* at $x$ if we can find an open set containing $x$ on which $f$ is constant. For $x_0 \in \mathcal{X}$, let $LC(x_0) = \{x : f(x) = f(x_0), f$ is locally constant at $x\}$. In [6], we provide the proof of the following result.

THEOREM 1. *Suppose that:*

(i) $\mathcal{X}$ *is a contented subset of $R^k$ with positive content;*

(ii) $B_n(x)$ *is a rectangle containing $x$ with $\mu_k(B_n(x))$ finite and constant in $x$ and such that $B_n(x)$ shrinks nicely to $x$ as $n \to \infty$;*

(iii) $\{B_n(x) : x \in R^k\}$ *forms a partition of $R^k$ with $\{B_{n+1}(x) : x \in R^k\}$ a subpartition of $\{B_n(x) : x \in R^k\}$ and $\sup_{x \in R^k} \operatorname{diam}(B_n(x)) \to 0$ as $n \to \infty$;*

(iv) $f$ *is a continuous density function on $\mathcal{X}$ with $f^{-1}A$ contented for any interval $A$ and such that $LC(x_0)$ is contented with $\mu_k(f^{-1}f(x_0) \cap LC(x_0)^c) = 0$.*

*Then (4) converges to $P(f(x) \le f(x_0))$ as $n \to \infty$.*

Theorem 1 gives conditions under which the appropriate discrete $P$-value, in the sense that we will always be measuring $x$ to some finite accuracy, is indeed approximated by the continuous version given by $P(f(x) \le f(x_0))$. Although the result will hold under weaker conditions, for example, we could allow $f$ to be piecewise smooth and allow for more general sets than rectangles in $\{B_n(x) : x \in R^k\}$, the conditions specified seem to apply in typical applications. For a specific example where $P(f(x) \leq f(x_0))$ fails to provide an approximation to (4) and where $f$ is a continuous density, see [6]. Basically, the approximation can fail when $f$ is too oscillatory so that (iv) does not hold. The example indicates, however, that this is more of a mathematical pathology than something we would encounter in a typical application. Condition (iv) can be substantially weakened if $P(f(x) = c) = 0$ for every $c$. In general, the distribution of $f(x)$ can have a discrete component, but our conditions imply that this can only arise by $f$ being constant on sets where it is locally constant. Also, Theorem 1 requires that the accuracy of the discretization is effectively the same across the sample space. In certain situations, we may want to allow this accuracy to vary across $\mathcal{X}$ so that we could obtain a different approximation to (4), but we do not pursue this issue further here.

We now consider basing the $P$-value on a general discrepancy statistic $T$. The question then is: given the discretization on $\mathcal{X}$ as determined by the measurement process, how should we take this into account? For, even if $T$ is 1–1, it will give rise to volume distortions and we do not want these to affect our $P$-value.

Suppose, first, that $\mathcal{X}$ and $\mathcal{T}$ are open subsets of $R^k$ and that $T$ is 1–1 and smooth. A partition element $B_n(x) \subset \mathcal{X}$ with measure $\mu_k(B_n(x))$ is then transformed into $TB_n(x)$ with measure $\mu_k(TB_n(x)) = \mu_k(B_n(x))J_T^{-1}(x')$ for some $x' \in B_n(x)$, while the density of the transformed response with respect to $\mu_k$ is $f_T(t) = f(T^{-1}(t))J_T(T^{-1}(t))$. Accordingly, we cannot use the $P$-value $P_T(f_T(t) \leq f_T(t_0))$ to assess whether or not $t_0$ or, equivalently, $x_0 = T^{-1}(t_0)$ is surprising since the density $f_T(t)$ depends on volume distortions and the sets $TB_n(x)$ are no longer necessarily of equal volume. There is clearly an easy fix for this, however, as we simply correct for this volume distortion and then $P_T(f_T(t)/J_T(T^{-1}(t)) \leq f_T(t_0)/J_T(T^{-1}(t))) = P(f(x) \leq f(x_0))$. With this refinement, (3) becomes invariant under 1–1, smooth transformations of the response $x$, that is, we retain as part of the problem prescription how the continuous probability model is approximating an essentially discrete response.

In general, however, $T$ will not be 1–1. Suppose, then, that $\mathcal{X}$ is an open subset of $R^k$ and $\mathcal{T}$ is an open subset of $R^l$, where $l \leq k$. Let $f_T$ denote the density of $P_T$ with respect to $\mu_l$ and suppose that this is continuous. Suppose that $T$ is sufficiently smooth so that for each $t \in \mathcal{T}$, the set $T^{-1}\{t\}$ is a smooth manifold with volume measure on $T^{-1}\{t\}$ denoted by $\nu_t$. As a simple example of this, suppose that $T(x) = x_1^2 + \cdots + x_k^2$ so that $T^{-1}\{t\}$ is a sphere of radius $t^{1/2}$, centered at the origin and $\nu_t$ is surface area measure. If $T$ is 1–1, then $T^{-1}\{t\}$ is a 0-dimensional manifold and $\nu_t$ is counting measure.

Results in [16] show that, in general,

$$(5) \qquad f_T(t) = \int_{T^{-1}\{t\}} f(x)|\det(dT(x) \circ dT'(x))|^{-1/2} \nu_t(dx),$$

where $dT$ is the differential of $T$. Formula (5) directly shows how $f_T$ is affected by volume distortions. For, at $x \in T^{-1}\{t\}$, the contribution to the density value $f_T(t)$ is distorted by the factor $J_T(x) = |\det(dT(x) \circ dT'(x))|^{-1/2}$. Accordingly, just as we do in the 1–1 case, we adjust the integrand in (5) by dividing by the factor $J_T(x)$ to obtain

$$f_T^*(t) = \int_{T^{-1}\{t\}} f(x)\nu_t(dx)$$

as the appropriate density to use. A simple example of this occurs when $T$ is $k$-to-one, so $T^{-1}\{t\} = \{x_1(t), \ldots, x_k(t)\}$ for each $t$. Then $\nu_t$ is counting measure and, by (5), $f_T(t) = \sum_{i=1}^{k} f(x_i(t)) J_T(x_i(t))$. In this case, we have that the corrected density is $f_T^*(t) = \sum_{i=1}^{k} f(x_i(t))$ and this equals $f(T^{-1}(t))$ when $k = 1$. In general, we see that $f_T^*$ is the density of $P_T$ with respect to the measure $(f_T(t)/f_T^*(t))\mu_l(dt)$ and the ratio $f_T(t)/f_T^*(t)$ measures the effect of the volume distortion induced by $T$ on the density $f_T$.

We then compute the $P$-value

$$(6) \qquad P_T(f_T^*(t) \le f_T^*(t_0)) = P(f_T^*(T(t)) \le f_T^*(T(x_0)))$$

to assess whether or not $t_0 = T(x_0)$ is a surprising value from $P_T$. This $P$-value depends only on the density assignment $f$ on the original response space, which is determined by how we are approximating an essentially discrete response, and the preimage sets of $T$.

We have the following simple, but significant, result for (6).

THEOREM 2. *If $\mathcal{X}$ is an open subset of $R^k$, $T : \mathcal{X} \to \mathcal{T}$ is onto with $\mathcal{T} \subset R^l$ open and $T$ is sufficiently smooth, then the $P$-value given by (6) is invariant under 1–1 smooth transformations on $\mathcal{T}$.*

PROOF. Suppose that $W$ is a 1–1, smooth transformation defined on $\mathcal{T}$ and that $w = W(t)$. Then $(W \circ T)^{-1}\{w\} = T^{-1}\{t\}$ and $f_{W \circ T}^*(w) = \int_{T^{-1}\{t\}} f(x) \times \nu_t(dx) = f_T^*(t)$. □

We now consider some examples and note that these support (6) as the appropriate definition of an invariant $P$-value.

EXAMPLE 1 ($P_T$ is discrete). First, suppose that $P$ is discrete on countable $\mathcal{X}$. Then $\nu_t$ is counting measure on $T^{-1}\{t\}$, $J_T(x) \equiv 1$ and, therefore, $f_T^*(t) = \int_{T^{-1}\{t\}} f(x)\nu_t(dx) = \sum_{x \in T^{-1}\{t\}} P(X = x) = p_T(t)$ is the probability function of $T$. Hence, (6) equals (2). Note that $dT$ is just the identity, so there is no volume distortion. If $P$ is continuous, then, for those $t$ with $p_T(t) > 0$, we have that $\nu_t$ is $\mu_k$ restricted to $T^{-1}\{t\}$. Accordingly, $f_T^*(t) = \int_{T^{-1}\{t\}} f(x)\nu_t(dx) = p_T(t)$ and, again, (6) equals (2).

EXAMPLE 2 ($J_T(x)$ is constant for $x \in T^{-1}\{t\}$). Note that $J_T(x)$ is constant for all $x$ whenever $T(x)$ is an affine transformation. Hence, we could have $T(x) = a + Bx$ for some $a \in R^l$ and $B \in R^{l \times k}$ when $\mathcal{X} \subset R^k$. Also, note that if $x \in R^n$ and $T(x)$ is the order statistic, then $J_T(x)$ is constant for all $x$. It is then clear that (6) equals $P_T(f_T(t) \le f_T(t_0))$. For example, when $T(x) = \bar{x}$, we simply use the density of $\bar{x}$ to compute the $P$-value so that when $P$ is the $N(0,1)$ distribution, (6) is $2(1 - \Phi(\bar{x}))$. As another example, suppose that $T$ is projection on the $i$th coordinate, so $J_T(x) \equiv 1$. Then $T^{-1}\{t\}$ is the set of points in $\mathcal{X}$ with $i$th coordinate equal to $t$, $\nu_t$ is Euclidean volume on this set and $f_T^*(t)$ is the marginal density of the $i$th coordinate. This generalizes to arbitrary coordinate projections.

The volume distortion can be constant in $T^{-1}\{t\}$ but vary with $t$. Putting $J_T^*(t) = J_T(x)$ for $x \in T^{-1}\{t\}$, from (5), we have that $f_T(t) = f_T^*(t)J_T^*(t)$, so (6) can be computed as $P_T(f_T(t)/J_T^*(t) \le f_T(t_0)/J_T^*(t_0))$. This permits us to avoid the integration involved in calculating $f_T^*(t)$ when we know the distribution of $T$ and can compute $J_T(x)$ easily.

For example, suppose that $T(x) = x'x$. Then $T^{-1}\{t\}$ is a sphere of dimension $(k-1)$ in $R^k$. Now, $dT(x) = 2(x_1 \cdots x_k)$, so $dT(x) \circ dT'(x) = 4x'x = 4t$ and $J_T(x) = 1/2t^{1/2}$ is constant for $x \in T^{-1}\{t\}$ for every $t$. Note that the adjustment factor involves multiplying $f_T(t)$ by $2t^{1/2}$ and this is precisely the distortion caused by the "quadratic" part of the transformation. The appropriate $P$-value is $P_T(f_T(t)t^{1/2} \le f_T(t_0)t_0^{1/2})$. We see that in this case, we must modify the usual density that we work with. As a particular case, suppose that $x \sim N_k(0, I)$. Then $T(x) \sim \chi^2(k)$ with density $f_T(t) = \Gamma^{-1}(k)2^{-k/2}t^{(k/2)-1}e^{-t/2}$. Therefore, the invariant $P$-value is given by $P_T(t^{(k-1)/2}e^{-t/2} \le t_0^{(k-1)/2}e^{-t_0/2})$ and only when $k = 1$ is this equivalent to $P_T(t \ge t_0)$. In contrast, when we directly observe $T \sim \chi^2(k)$, in the sense that it is a measured variable, and we discretize using equal length intervals, the relevant $P$-value is $P_T(t^{(k/2)-1}e^{-t/2} \le t_0^{(k/2)-1}e^{-t_0/2})$. As was just shown, when we take into account that $T$ arises as a transformation of a measured variable, the $P$-value changes. Further, both of these $P$-values are two-sided when $k > 1$.

As we will see, Example 3 is a situation where $J_T(x)$ varies with $x \in T^{-1}\{t\}$.

In Section 3, we discuss a situation that involves comparing the observed $x_0$ with the conditional distribution of $x$ given that $W(x) = W(x_0) = w_0$ for some smooth transformation $W$. In this case, the conditional density of $x$, with respect to volume measure on $W^{-1}\{w_0\}$, is $f(x)J_W(x)/f_W(w_0)$ and it is clear that the volume distortion at $x$, induced by the conditioning, is given by $J_W(x)$. Accordingly, the relevant $P$-value, based on the full data, is given by $P(f(x_0)/f_W(w_0) \leq f(x)/f_W(w_0)|W(x) = w_0) = P(f(x) \leq f(x_0)|W(x) = w_0)$. If we have a transformation $T$ of $x$, then the relevant $P$-value is as follows.

LEMMA 3.  *Suppose that* $\mathcal{X}, \mathcal{W}$ *and* $\mathcal{T}$ *are manifolds with volume measures* $\mu_{\mathcal{X}}, \mu_{\mathcal{W}}$ *and* $\mu_{\mathcal{T}}$, *respectively, and* $W: \mathcal{X} \to \mathcal{W}, T: \mathcal{X} \to \mathcal{T}$ *are onto and smooth. Let* $\nu_{T,W,t,w}$ *denote volume measure on* $T^{-1}\{t\} \cap W^{-1}\{w\}$. *The relevant conditional P-value based on* $T$, *given* $W(x) = w_0$, *is then*

$$(7) \qquad P_T(f^*_{T,W}(t|w_0) \leq f^*_{T,W}(t_0|w_0)|W(x) = w_0),$$

*where* $t_0 = T(x_0)$ *and* $f^*_{T,W}(t|w_0) = \int_{T^{-1}\{t\} \cap W^{-1}\{w_0\}} f(x)\nu_{T,W,t,w_0}(dx)$.

PROOF.  The conditional density of $T$ given $W = w$ is given by $f_{T,W}(t|w) = \int_{T^{-1}\{t\} \cap W^{-1}\{w\}} (f(x)/f_W(w))J_{(T,W)}(x)\nu_{T,W,t,w}(dx)$. So, the volume distortion induced by the transformations is $J_{(T,W)}(x)$ and the result follows.  □

We will also need a technical result concerning the composition of mappings.

LEMMA 4.  *Suppose that* $\mathcal{X}, \mathcal{U}$ *and* $\mathcal{T}$ *are manifolds with volume measures* $\mu_{\mathcal{X}}, \mu_{\mathcal{U}}$ *and* $\mu_{\mathcal{T}}$, *respectively, and* $U: \mathcal{X} \to \mathcal{U}, T: \mathcal{U} \to \mathcal{T}$ *are onto, smooth mappings. Then*

$$f^*_{T \circ U}(t) = \int_{T^{-1}\{t\}} J_T(u) \int_{U^{-1}\{u\}} f(x) J^{-1}_{T \circ U}(x) J_U(x)\nu_{U,u}(dx)\nu_{T,t}(du),$$

*where* $\nu_{U,u}$ *and* $\nu_{T,t}$ *are the volume measures on* $U^{-1}\{u\}$ *and* $T^{-1}\{t\}$, *respectively.*

PROOF.  Suppose that $g: \mathcal{X} \to R^1$ is nonnegative, $\int_A g(x)\mu_{\mathcal{X}}(dx)$ is finite for compact $A$ and $B \subset \mathcal{T}$ is open. By the measure decomposition theorem (see [16], Theorem 15.1) applied to $g(x)\mu_{\mathcal{X}}(dx)$ and $T \circ U$, we have that $\int_{\mathcal{X}} I_B(T(U(x)))g(x)\mu_{\mathcal{X}}(dx) = \int_B \int_{(T \circ U)^{-1}\{t\}} g(x)J_{T \circ U}(x)\nu_{T \circ U,t}(dx)\mu_{\mathcal{T}}(dt)$. Apply the measure decomposition theorem to $I_B(T(U((x)))\mu_{\mathcal{X}}(dx)$ and $U$, and then to $\int_{U^{-1}\{u\}} I_B(T(U(x)))g(x)J_U(x)\nu_{U,u}(dx)\mu_{\mathcal{U}}(du)$ and $T$, to obtain $\int_{\mathcal{X}} I_B(T(U(x)))g(x)\mu_{\mathcal{X}}(dx) \qquad = \qquad \int_{\mathcal{U}} \int_{U^{-1}\{u\}} I_B(T(U(x)))g(x)J_U(x) \qquad \times$ $\nu_{U,u}(dx)\mu_{\mathcal{U}}(du) = \int_B \int_{T^{-1}\{t\}} \int_{U^{-1}\{u\}} g(x) J_U(x) \nu_{U,u}(dx) J_T(u) \nu_{T,t}(du) \mu_{\mathcal{T}}(dt)$. Then $\int_{(T \circ U)^{-1}\{t\}} g(x) J_{T \circ U}(x) \nu_{T \circ U,t}(dx) \mu_{\mathcal{T}}(dt) = \int_{T^{-1}\{t\}} \int_{U^{-1}\{u\}} g(x) J_U(x) \times \nu_{U,u}(dx) J_T(u) \nu_{T,t}(du) \mu_{\mathcal{T}}(dt)$, and setting $g(x) = f(x) J_{T \circ U}^{-1}(x)$ establishes the result. $\square$

**3. Model checking.** Suppose that we have a statistical model $\{P_\theta : \theta \in \Theta\}$, where $P_\theta$ is a probability measure on $\mathcal{X}$ with density $f_\theta$ with respect to support measure $\mu_k$. We now investigate $P$-values for checking the model in light of the observed $x_0$.

If $W : \mathcal{X} \to \mathcal{W}$ is a minimal sufficient statistic, then the conditional distribution of the data given $W$ is independent of $\theta$ and is denoted $P(\cdot|W(x) = w_0)$. To check the model, we can then compare $x_0$ to $P(\cdot|W(x) = w_0)$ to see if the observed data is surprising. By the converse of the factorization theorem, we have that $f_\theta(x) = g_\theta(W(x)) h(x)$. Lemma 3 and (7) give an invariant $P$-value that assesses $f_\theta$ for each $\theta$. We have the following result.

THEOREM 5. *For the statistical model $\{P_\theta : \theta \in \Theta\}$ with minimal sufficient statistic $W$, the $P$-value (7), associated with discrepancy statistic $T$, equals*

$$(8) \qquad P_T(h_{T,W}(t|w_0) \le h_{T,W}(t_0|w_0)|W(x) = w_0),$$

*where $h_{T,W}(t|w_0) = \int_{T^{-1}\{t\} \cap W^{-1}\{w_0\}} h(x) \nu_t(dx)$. That is, it is independent of $\theta$ and (8) is independent of the choice of $h$.*

PROOF. In the continuous case, we assume that each density is continuous at any observed $x_0$ and restrict our attention to those $x_0$ for which $f_\theta(x_0) > 0$. If $f_\theta(x_0) > 0$, then $g_\theta(W(x_0)) > 0$ and $g_\theta(W(x)) = g_\theta(W(x_0))$ for the event $W(x) = t_0 = W(x_0)$. We have that (7) equals

$$P_T \left( \begin{array}{c} \int_{T^{-1}\{t\} \cap W^{-1}\{w_0\}} f_\theta(x) \nu_t(dx) \\ \le \int_{T^{-1}\{t_0\} \cap W^{-1}\{w_0\}} f_\theta(x) \nu_{t_0}(dx) \end{array} \middle| W(x) = w_0 \right)$$

$$= P(h_{T,W}(t|w_0) \le h_{T,W}(t_0|w_0)|W(x) = w_0).$$

Further, if $g_\theta(W(x)) h(x) = g'_\theta(W(x)) h'(x)$, then

$$P_T(h_{T,W}(t|w_0) \le h_{T,W}(t_0|w_0)|W(x) = w_0)$$

$$= P_T \left( \begin{array}{c} \int_{T^{-1}\{t\} \cap W^{-1}\{w_0\}} g_\theta(W(x)) h(x) \nu_t(dx) \\ \le \int_{T^{-1}\{t_0\} \cap W^{-1}\{w_0\}} g_\theta(W(x)) h(x) \nu_{t_0}(dx) \end{array} \middle| W(x) = w_0 \right)$$

$$= P_T \left( \left. \begin{array}{l} \displaystyle\int_{T^{-1}\{t\} \cap W^{-1}\{w_0\}} g'_\theta(W(x))h'(x)\nu_t(dx) \\ \displaystyle\leq \int_{T^{-1}\{t_0\} \cap W^{-1}\{w_0\}} g'_\theta(W(x))h'(x)\nu_{t_0}(dx) \end{array} \right| W(x) = w_0 \right)$$

$$= P_T(h'_{T,W}(t|w_0) \leq h'_{T,W}(t_0|w_0) | W(x) = w_0)$$

and we are done. $\square$

We now consider an application of this result.

EXAMPLE 3 (Checking a normal model using the Jarque–Bera test statistic). Suppose that $x = (x_1, \ldots, x_n)$ is a sample of $n$ from the $N(\mu, \sigma^2)$ distribution with $\mu \in R^1$ and $\sigma^2 > 0$ unknown. Then $W(x) = (\bar{x}, r)$, where $r = \|x - \bar{x}1_n\|$, is minimal sufficient. Putting $d(x) = (x - \bar{x}1_n)/r$, we can write $x = \bar{x}1_n + rd$ and note that $\bar{x}, r$ and $d$ are statistically independent with $d$ uniformly distributed on $S^{n-1} \cap L^\perp\{1_n\}$. In this case, $h$ is constant (so we can take it to be 1) and $W^{-1}\{(\bar{x}_0, r_0)\}$ is the $(n-2)$-dimensional sphere $\bar{x}_0 1_n + r_0(S^{n-1} \cap L^\perp\{1_n\})$.

It is natural here to consider functions of $d$ as discrepancy statistics for checking the model. If $T$ is a real-valued function of $d$, then $h_{T,W}(t|w_0)$ is the volume of the $(n-3)$-dimensional submanifold of $\bar{x}_0 1_n + r_0(S^{n-1} \cap L^\perp\{1_n\})$ given by $(T \circ d)^{-1}\{t\} \cap W^{-1}\{(\bar{x}_0, r_0)\}$. Alternatively, from the proof of Theorem 5, we can compute the invariant $P$-value by assuming $(\mu, \sigma) = (0, 1)$, letting $f$ denote the density of a sample of $n$ from the $N(0, 1)$ distribution and computing $P_{(0,1)}(f^*_{T \circ d}(T(d(x))) \leq f^*_{T \circ d}(T(d(x_0)))|(\bar{x}_0, r_0)) = P_d(f^*_{T \circ d}(T(d)) \leq f^*_{T \circ d}(T(d_0)))$, where $f^*_{T \circ d}(t) = \int_{(T \circ d)^{-1}\{t\}} f(x)\nu_t(dx)$ and $d$ is uniformly distributed on $S^{n-1} \cap L^\perp\{1_n\}$.

As an illustration, we consider a commonly used discrepancy statistic for checking normality, namely, the Jarque–Bera statistic $T = n(nT_3^2/6 + (nT_4 - 3)^2/24)$, where $T_p \circ d = \sum_{i=1}^n d_i^p$. The statistic $T$ is an attempt to create an omnibus test by combining the skewness and kurtosis statistics. The form is based on the asymptotic normality of $(T_3, T_4)$ as this implies that $T$ is asymptotically $\chi^2(2)$.

The volume distortion induced by $T_p(d)$ can be computed explicitly as $J_{T_p \circ d}(x) = (r/p)(T_{2p-2}(d(x)) - T_{p-1}^2(d(x))/n - T_p^2(d(x)))^{-1/2}$. From this, we obtain that $J_{T \circ d}(x)$ for the Jarque–Bera statistic $T$ satisfies $J_{T \circ d}(x)^{-2} = (n^4/r^2)[(nT_4(d(x))/3 - 1)^2 T_6(d(x)) + 2(nT_4(d(x))/3 - 1)T_3(d(x))T_5(d(x)) - (T_3^2(d(x)) + nT_4^2(d(x))/3 - T_4(d(x)))^2 + T_3^2(d(x))T_4(d(x)) - nT_3^2(d(x)) \times T_4^2(d(x))/9]$. We see that $J_T(x)$ is not a function of $T \circ d$ and so the volume distortion is not constant within $T^{-1}\{t\}$.

In Figure 1, we have plotted the densities $f_T$ and the invariant $P$-values based on the Jarque–Bera test statistic for several sample sizes $n$. The densities are quite irregular for small sample sizes and skewed. Note that, while the formula for the Jarque–Bera statistic is simple, it is, in fact, a degree 8 polynomial in the $d_i$. The volume distortion reflects this complexity and this is seen to have an effect, even for very large sample sizes. We estimated $f_T(t)$ and $f_T^*(t) = f_T(t)E(J_T^{-1}(X)|T=t)$ via simulation using kernel density estimation methods.

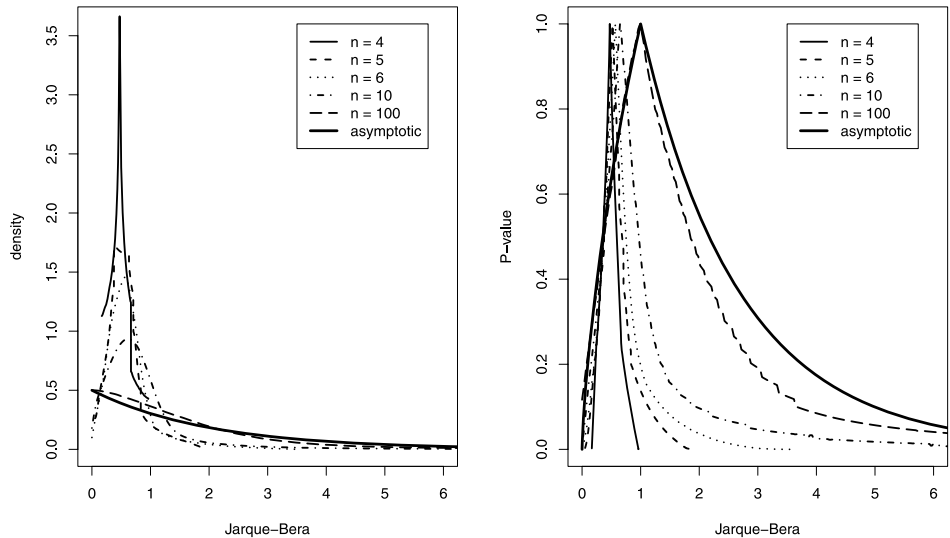

Fig. 1. *Densities and invariant P-values for Jarque–Bera test for various sample sizes n when sampling from normal.*

There is no reason to suppose, however, that the Jarque–Bera statistic represents the best way to combine $T_3$ and $T_4$. It seems more natural to set $T = (T_3, T_4)$ and then compute (6), as this takes account of the dependence between $T_3$ and $T_4$. As illustrated in Figure 2, the $P$-values obtained via (6) and (3) with this $T$ are very similar, indicating that volume distortion is playing a very small role, even for small sample sizes. Also, these $P$-values are much more critical of the normality assumption than the Jarque–Bera test and we note that these are based on the precise form of the joint distribution of $T_3$ and $T_4$. Finally, as $n$ increases, these $P$-values converge to the same $P$-value as given by (1) using the Jarque–Bera statistic. Overall, this approach to combining $T_3$ and $T_4$ seems to have considerable advantages over the Jarque–Bera test.

A similar analysis can be carried out for other discrepancy statistics. For example, the Shapiro–Wilk test has $J_T(x)$ constant for $x \in T^{-1}\{t\}$. Correcting for volume distortion results in a small change in the usual $P$-values for small $n$ and this effect disappears as $n$ grows. Overall, the Shapiro–Wilk test is much better behaved than the Jarque–Bera test. Note that the Shapiro–Wilk test statistic is a quadratic function of the $d_i$.

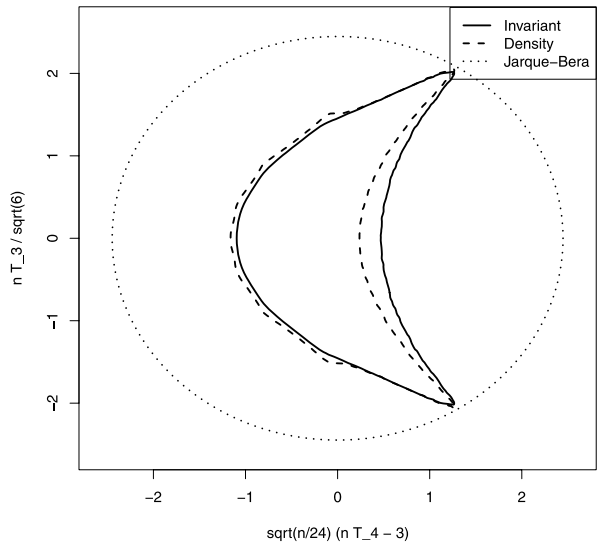

Fig. 2. *Contour of $(T_3, T_4)$ giving the 0.05 point for tests based on (6), (3) and the Jarque–Bera test, when $n = 10$.*

When $U$ is an ancillary statistic, a function of $U$ can be used to assess the model. If we consider the transformation $T \circ U$, then we must evaluate $\int_{(T \circ U)^{-1}\{t\}} f_\theta(x)\nu_t(dx)$ and, in general, there is no reason to suppose that this is independent of $\theta$. If the distribution of $T \circ U$ is discrete, however, then $\int_{(T \circ U)^{-1}\{t\}} f_\theta(x)\nu_t(dx)$ is the probability function of $T \circ U$ and, as such, is independent of $\theta$. Theorem 6 will show that the $P$-value based on $\int_{(T \circ U)^{-1}\{t\}} f_\theta(x) \times \nu_t(dx)$ is independent of $\theta$ for a very broad class of ancillaries.

Consider the following example which will serve as an archetype for a common situation where ancillaries arise.

EXAMPLE 4 (Location-scale models). Suppose that we have $x \in R^n$ and that the model is $x = \mu 1_n + \sigma z$, where $z$ is distributed with density $f$ with respect to volume measure on $R^n$ and $\mu \in R^1, \sigma > 0$ are unknown. Then $x$ has density $f_{\mu,\sigma}(x) = \sigma^{-n}f((x - \mu 1_n)/\sigma)$. We take the parameter space to be $\Theta = \{(\mu, \sigma) : \mu \in R^1, \sigma > 0\}$ and note that we have a group product defined on $\Theta$ via $(\mu_1, \sigma_1)(\mu_2, \sigma_2) = (\mu_1 + \sigma_1\mu_2, \sigma_1\sigma_2)$. This group acts on $R^n$ via $(\mu, \sigma)x = \mu 1_n + \sigma x$. A maximal invariant is then given by $U(x) = (x - x_m^* 1_n)/(x_{q3}^* - x_{q1}^*)$, where $x_m^*$ is the sample median, that is, the middle sample value, when $n$ is odd and the average of the two middle sample values when $n$ is even, and $x_{qi}^*$ is the $i$th quartile. Since $U(x)$ is invariant, it is ancillary. Note that $U^{-1}\{u\} = \{x : x = a 1_n + cu$ for some $(a, c) \in \Theta\} = \Theta u$, that is, $U^{-1}\{u\}$ is an orbit of the group action. Clearly, this orbit is half of a two-dimensional plane in $R^n$ and so the volume measure $\nu_u$ is just area.

If we wish to base our checking on $U$ itself, then we must evaluate

$$f^*_{\mu,\sigma,U}(d_0) = \int_{U^{-1}\{u\}} f_{\mu,\sigma}(x)\nu_u(dx) = \int_0^\infty \int_{-\infty}^\infty f_{\mu,\sigma}(a1_n + cu)\sqrt{n}\,da\,dc$$

$$= \int_0^\infty \int_{-\infty}^\infty \sigma^{-n} f\left(\frac{a-\mu}{\sigma}1_n + \frac{c}{\sigma}u\right)\sqrt{n}\,da\,dc$$

$$= \sigma^{-(n-2)} \int_0^\infty \int_{-\infty}^\infty f(a1_n + cu)\sqrt{n}\,da\,dc.$$

Accordingly, the $P$-value for model checking is given by

$$P_U(f^*_{\mu,\sigma,U}(u) \le f^*_{\mu,\sigma,U}(u_0))$$

$$= P_U\left(\int_0^\infty \int_{-\infty}^\infty f(a1_n + cu)\,da\,dc \le \int_0^\infty \int_{-\infty}^\infty f(a1_n + cu_0)\,da\,dc\right)$$

and, as this is independent of the model parameter, we have a valid $P$-value for checking the model. If, instead, we use a function $T(U)$, then an application of Lemma 4 shows that $f^*_{\mu,\sigma,T\circ U}$ is independent of $(\mu,\sigma)$, by the same argument, as the Jacobian factors do not depend on the parameter. Note that if $f$ is the $N(0,1)$ density, then basing model checking on the ancillary $d$ or on the conditional distribution of the data given a minimal sufficient statistic produce the same results.

More generally, suppose we have a group model $\{f_g : g \in G\}$, where $G$ is a group with a smooth product acting freely and smoothly on $\mathcal{X}$ and $f_g(x) = f(g^{-1}x)J_g(g^{-1}x)$ for some fixed density $f$. Now, suppose that $[\cdot] : \mathcal{X} \to G$ is smooth and satisfies $[gx] = g[x]$ so that $U(x) = [x]^{-1}x$ is a maximal invariant and is thus ancillary. Hence, $u = U(x) \in \mathcal{X}, x = [x]U(x)$ and $U^{-1}\{u\}$ is the orbit $\{gu : g \in G\}$. Now, if $\nu_u^*$ denotes volume measure on $G$, we have that $\nu_u = K(u)\nu_G^*$ for some positive function $K$. Let $z = g^{-1}x$ so that $[z] = g^{-1}[x]$ and let $J_g^*([z])$ denote the Jacobian of the transformation $[z] \to [x]$. We then have $f^*_{gU}(u) = \int_{U^{-1}\{u\}} f_g(x)\nu_u(dx) = \int_{\{gu : g \in G\}} f_g([x]u)\nu_u(dx) = K(u)\int_G f(g^{-1}[x]u)J_g(g^{-1}[x]u)\nu_G^*(d[x]) = K(u)\int_G f([z]u)J_g(u)J_g^*([z])\nu_G^*(d[z])$. Now, if we can write $J_g(u)J_g^*(z) = L(u)m(g)$ for some positive functions $L$ and $m$, then we have that the invariant $P$-value $P_U(f^*_{gU}(u) \le f^*_{gU}(u_0))$ is indeed independent of $g$. Further, by Lemma 4, this will also hold for $T \circ U$. In Example 4, $J_{(\mu,\sigma)}(u) = \sigma^{-n}$ and, with $[x] = [x^*_m, x^*_{q3} - x^*_{q1}]$, we have $J_g^*([z]) = \sigma^2$ and this condition is satisfied. More generally, this condition is satisfied in a wide range of group models, such as those discussed in [7]. Accordingly, the following result is broadly applicable.

THEOREM 6. *Suppose that $\{f_g : g \in G\}$ is a family of densities with respect to volume measure $\mu_{\mathcal{X}}$ on $\mathcal{X}$, where $G$ is a group with a smooth product and a smooth action defined on $\mathcal{X}$, and $f_g(x) = f(g^{-1}x)J_g(g^{-1}x)$. Further, suppose that there exists a smooth $[\cdot] : \mathcal{X} \to G$ satisfying $[gx] = g[x]$ and let $J_g^*([z])$ denote the Jacobian of the transformation $[z] \to [x]$, where $x = gz$. If there exist positive functions $L$ and $m$ such that $J_g(u)J_g^*(z) = L(u)m(g)$, then we have that the P-value (8) based on the ancillary $T \circ U$, with $U(x) = [x]^{-1}x$ and $T$ smooth, is independent of the model parameter and is thus a valid check on the model.*

## 4. Conclusions.
The use of $P$-values is a somewhat controversial topic. When we are concerned with model checking, however, their use seems unavoidable. At least one problem with $P$-values is the ambiguity as to how they should be defined when we have data $x_0$, a single probability measure $P$ from which the data was supposedly generated and a discrepancy statistic $T$. While (1) has some appeal, the rationale for this is dependent on the form of the distribution $P_T$, namely, the right tail being the only region where surprising values of $T(x_0)$ may occur, and this clearly does not hold for an arbitrary $T$. For the discrete case, (2) seems like a much more appealing definition. A difficulty with (2) arises when we try to generalize this to absolutely continuous contexts since the simple analog of (2) is not invariant under 1–1, smooth transformations of $T$.

We have argued that an appropriate, generalized definition of an invariant $P$-value can be obtained from (2). For this, we must take into account that a continuous model is essentially an approximation to a true discrete model as we always measure responses to some finite accuracy. Further, we must not allow volume distortions induced by a discrepancy statistic $T$ to have any effect on our inferences. This leads us to the invariant $P$-value given by (6). The definition is seen to be intuitively reasonable and to perform well in a variety of examples. For other model checking problems, where our approach to $P$-values seems applicable, see [1, 2, 8] and [11].

**Acknowledgments.** We wish to thank the reviewers for a number of suggestions that lead to improvements in this paper.

DEPARTMENT OF STATISTICS
UNIVERSITY OF TORONTO
TORONTO, ON M5S 3G3
CANADA
E-MAIL: mevans@utstat.utoronto.ca
           gunho@utstat.utoronto.ca